\documentstyle[12pt, amstex, amsfonts]{amsart}  

\title{Minimal Sequences and the Kadison-Singer Problem}

\author{Wayne Lawton}

\address{Department of Mathematics\\
    National University of Singapore \\
    10 Lower Kent Ridge Road\\
    Singapore 119076
    email: matwml(AT)nus.edu.sg}
\newtheorem{thm}{Theorem}[section]
\newtheorem{lem}{Lemma}[section]

\newtheorem{cor}{Corollary}[section]

\newtheorem{defin}{Definition}[section]

\def\L2{L^2([ 0,2\pi))}
\def\L2R{L^2(\bold{R})}

\def\cA{{\cal A }}
\def\cB{{\cal B }}

\def\cL{{\cal L }}
\def\cM{{\cal M }}

\def\cP{{\cal P }}

\def\CC{{\Bbb C}}
\def\NN{{\Bbb N}}
\def\QQ{{\Bbb Q}}
\def\RR{{\Bbb R}}
\def\TT{{\Bbb T}}

\def\ZZ{{\Bbb Z}}

\begin{document} 

\begin{abstract}
The Kadison-Singer problem asks: does every pure state on the
C$^*$-algebra $\ell^{\infty}(\ZZ)$ admit a unique extension to the
C$^*$-algebra $\cB(\ell^2(\ZZ))$? A yes answer is equivalent to
several open conjectures including Feichtinger's: every bounded
frame is a finite union of Riesz sequences. We prove that for
measurable $S \subset \TT,$
$\{ \chi_{_S} \, e^{2\pi i k t} \}_{_{k\in \ZZ}}$
is a finite union of Riesz sequences in $L^2(\TT)$ if and only if
there exists a nonempty $\Lambda \subset \ZZ$ such that
$\chi_{_\Lambda}$ is a minimal sequence and
$\{ \chi_{_S} \, e^{2\pi i k t} \}_{_{k \in \Lambda}}$
is a Riesz sequence. We also suggest some directions for future
research.
\end{abstract}

\maketitle

\section{Introduction}
\setcounter{equation}{0}
Recently there has been considerable interest in two deep problems
that arose from very separate areas of mathematics. The

\medskip

\noindent {\bf Kadison-Singer Problem (KSP):} Does every pure state
on the C$^*$-algebra $\ell^{\infty}(\ZZ)$ admit a unique extension to
the C$^*$-algebra $\cB(\ell^2(\ZZ))$?

\medskip

\noindent arose in the area of operator algebras and has remained
unsolved since 1959 \cite{KS59}. Pure states correspond to points in
a topological space, the Stone-\v{C}ech compactification $\beta
(\ZZ)$ of $\ZZ,$ whose construction requires the axiom of choice,
and recent work implicates the KSP with set-theoretic foundational
issues \cite{S08}. The

\medskip

\noindent {\bf Feichtinger Conjecture (FC):} Every bounded frame can
be written as a finite union of Riesz sequences.

\medskip

\noindent arose from Feichtinger's work in the area of signal
processing involving time-frequency analysis, \cite{GR00},
\cite{GR03}, (\cite{CH03},References) and has remained unsolved
since it was formally stated in the literature in 2005 (\cite{CA05},
Conjecture 1.1).

Casazza and Tremain proved (\cite{CA06b}, Theorem 4.2) that a yes
answer to the KSP is equivalent to the FC and Casazza, Fickus,
Tremain, and Weber explained many other equivalent conjectures in
\cite{CA06a}. In this paper we address the

\medskip

\noindent {\bf Feichtinger Conjecture for Exponentials (FCE):} For
every non-trivial measurable set $S \subset \TT,$ the sequence
$\{ \chi_{_S} \, e^{2\pi i k t} \}_{_{k\in \ZZ}}$
is a finite union of Riesz sequences.

\medskip

\noindent Although FC implies FCE, and FCE is easily shown to be
equivalent to FC for frames of translates, it is unknown if FCE
implies FC. Our intuition suggests that FCE is weaker than FC. Our
main result relates FCE to the area of {\bf Symbolic Dynamics:}
\begin{thm}
\label{thm_main}
For subsets $S \subset \TT$ and $\Lambda \subset \ZZ$ set
$B(S,\Lambda) := \{ \chi_{_S} \, e^{2\pi i k t} \}_{_{k\in \ZZ}}.$
For every nontrivial measurable $S \subset \TT$ the following
conditions are equivalent:
\begin{enumerate}
\item $B(S,\ZZ)$ is a finite union of Riesz sequences,
\item there exists a {\bf syndetic} subset $\Lambda \subset \ZZ$
such that $B(S,\Lambda)$ is a Riesz sequence,
\item there exists a nonempty subset $\Lambda \subset \ZZ$
such that $\chi_{_\Lambda}$ is a {\bf minimal sequence} and
$B(S,\Lambda)$ is a Riesz sequence.
\end{enumerate}
\end{thm}
The remainder of this section introduces notation, derives
preliminary results, and reviews selected known results. Section 2
derives Theorems (\ref{thm_main}) and (\ref{thm2}). Section 3
suggests some directions for further research.
$\NN = \{1,2,...\},$ $\ZZ,$ $\QQ,$ $\RR,$ $\CC$
are the natural, integer, rational, real, and complex numbers,
$\TT = \RR / \ZZ$ is the circle group, $\cL^{+}(\TT)$ is the set of Lebesque
measurable $S \subseteq \TT$ whose Haar measure $\mu(S) > 0,$ and
$F_n := \{0,1,...,n-1\}.$
For $Y \subset X,$ $X \backslash Y$ is the complement of
$Y$ in $X$ and $\chi_{_Y} \, : \, X \rightarrow \{0,1\}$ is the
characteristic function of $Y.$
For $S \in \cL^{+}(\TT)$ and $\Lambda \subset \ZZ,$ $P_{_S},
P_{_\Lambda}$ are orthogonal projections of $L^2(\TT)$ onto the
closed subspace $\chi_{_S} \, L^2(\TT),$ the closed subspace spanned
by the sequence $E(\Lambda) := \{ e^{2\pi i k \, t } \}_{_{k \in
\Lambda}},$ respectfully.
\begin{lem}
\label{lem1.1}
The following conditions are equivalent:
\begin{enumerate}
\item $\exists \, \epsilon_1 > 0$ such that
    $|| \, P_{_S} \, P_{_\Lambda} \, h \, || \geq
    \epsilon_1 \, || \, P_{_\Lambda} \, h \, ||, \ h \in L^2(\TT),$
\item $\exists \, \epsilon_2 > 0$ such that
    $|| \, P_{_S} \, h \, || + || \, P_{_{\ZZ \backslash \Lambda}} \, h \, || \geq
    \epsilon_2 \, || \, h \, ||, \ h \in L^2(\TT),$
\item $\exists \, \epsilon_3 > 0$ such that
    $|| \, P_{_{\ZZ \backslash \Lambda}} \, P_{_{\TT \backslash S}} \, h \, || \geq
    \epsilon_3 \, || P_{_{\TT \backslash S}} \, h \, ||,
    \ h \in L^2(\TT).$
\end{enumerate}
\end{lem}
{\begin{pf} Clearly (2) implies (1) and (3).
Let $h \in L^2(\TT).$ Then $h = h_1 \cos \theta + h_2 \sin \theta$ where $\theta \in [0,\frac{\pi}{2}],$
$h_1 \cos \theta = P_{_\Lambda} \, h,$ $h_2 \sin \theta = P_{_{\ZZ \backslash \Lambda}} \, h,$ and
$||h_1|| = ||h_2|| = ||h||.$
Hence (1) implies
$
||P_{_S}h|| + ||P_{_{\ZZ \backslash \Lambda}}\widehat h|| \geq
\left(\max\{ 0, \epsilon_1 \cos \theta - \sin \theta \} + \sin \theta \right) \, ||h||
$
so (2) holds with $\epsilon_2 = \epsilon_1 \, \left( \, 1 + \epsilon_{1}^2 \, \right)^{-1/2}.$
A similar argument shows that (3) implies (2).
\end{pf}
Christenson's book \cite{CH03} explains frames
and Riesz sequences. $B(S,\Lambda)$ is a bounded (below by $|S|$) frame in $P_S L^2(\TT).$ In Lemma
(\ref{lem1.1}) condition (1) holds iff $B(S,\Lambda)$ is a Riesz sequence and
condition (3) holds iff $R_{_{\TT \backslash S}}E(\ZZ \backslash \Lambda)$ is
a frame in $L^2(\TT \backslash S).$ Here
$R_{_{\TT \backslash S}} : L^2(\TT) \rightarrow L^2(\TT \backslash S)$
is the restriction operator.

\medskip

\noindent For $\Lambda \subset \ZZ$ we define lower and upper Beurling densities
$$D^{-}(\Lambda) = \lim_{k \rightarrow \infty} \min_{a \in \RR} \frac{|\Lambda \cap (a,a+k)|}{k}, \ \ \ \
D^{+}(\Lambda) = \lim_{k \rightarrow \infty} \max_{a \in \RR} \frac{|\Lambda \cap (a,a+k)|}{k},$$
lower and upper asymptotic densities
$$d^{-}(\Lambda) = \liminf_{k \rightarrow \infty} \frac{|\Lambda \cap (-k,k)|}{2k}, \ \ \ \
d^{+}(\Lambda) = \limsup_{k \rightarrow \infty} \frac{|\Lambda \cap (-k,k)|}{2k},$$
and if the cardinality $|\Lambda| \geq 2$ we define the separation
$$\Delta(\Lambda) := \min \{ \, |\lambda_2-\lambda_1| \, : \, \lambda_1, \lambda_2 \in \Lambda, \, \lambda_1 \neq
\lambda_2 \, \}.$$
The following result was inspired by Olevskii and Ulanovskii's paper
\cite{O08}:
\begin{cor}
\label{cor1.1}
If $B(S,\Lambda)$ is a Riesz set then $D^{+}(\Lambda) \leq \mu(S).$
\end{cor}
\begin{pf}
Since $E(\ZZ \backslash \Lambda)$ is a frame in $L^2(\TT \backslash
\Lambda)$ Landau's result (\cite{L67}, Theorem 3) implies
$D^{-}(\ZZ \backslash \Lambda) \geq \mu(\TT \backslash S).$
Therefore
$D^{+}(\Lambda) = 1 - D^{-1}(\ZZ \backslash \Lambda) \leq 1 -
\mu(\TT \backslash S) = \mu(S).$
\end{pf}
{\bf Result 1} Montgomery and Vaughan's result (\cite{MV74},
Corollary 2) implies that if $S$ contains an interval having length
$T > 1/\Delta(\Lambda)$ then condition (1) in Lemma (\ref{lem2.1})
holds with $\epsilon_1 = T - 1/\Delta(\Lambda)$ so $B(S,\Lambda)$ is
a Riesz sequence. It follows that if $B(S,\ZZ)$ does not satisfy FCE
then there exists a Cantor set $S_c \in \cL^{+}(\TT)$ such that $S_c
\subseteq S.$

\medskip

\noindent {\bf Result 2} Casazza, Christiansen, and Kalton showed
(\cite{CA01}, Theorem 2.2) that for $n \in \NN, m \in \ZZ,$
$B(S,n\ZZ + m)$ is a Riesz basis iff $S + \frac{1}{n} F_n = \TT$
a.e. This condition never holds if $S$ is a Cantor set.

\medskip

\noindent {\bf Result 3} The authors above also showed (\cite{CA01},
Theorem 2.4) that for $\Lambda \subseteq \NN,$ $B(S,\Lambda)$ is a
Riesz sequence iff $B(S,\Lambda)$ is a frame.

\medskip

\noindent {\bf Result 4} Bourgain and Tzafriri's restricted
invertibility result for matrices \cite{BT87} implies that for every $S
\in \cL^{+}(\TT)$ there exists $\Lambda \subseteq \ZZ$ such that
$d^{-}(\Lambda) > 0$ and $B(S,\Lambda)$ is a Riesz sequence.

\medskip

\noindent {\bf Result 5} Bourgain and Tzafriri's result
(\cite{BT91}, Theorem 4.1) implies that if $\chi_{_S}$ belongs to
the Besov space $W_{2,2}^{\tau}$ for some $\tau > 0$ then $B(S,\ZZ)$
satisfies FCE. Moreover, the proof of their result (\cite{BT91},
Corollary 4.2) shows that if $S$ is a Cantor set and $\TT \backslash
S$ is a union of disjoint open intervals $I_n, \, n \in \NN$
satisfying $\mu(I_n) \leq c 2^n$ for some $c > 0$ then $\chi_{_S}
\in W_{2,2}^{\tau}$ for all $\tau \in (0,1).$

\medskip

\noindent {\bf Result 6} Bownik and Speegle (\cite{BS06}, Theorem
4.16) used discrepancy theory to construct $S \in \cL^{+}(\TT)$ and
a class of $\Lambda \subset \ZZ$ such that $B(S,\Lambda)$ is not a
Riesz sequence and related their construction to Gower's results
about Szemeredi's theorem \cite{GO01}.

\medskip

\noindent {\bf Result 7} In November 2009 Spielman and N. Srivastava
gave an elementary constructive proof of Bourgain and Tzafriri's
restricted invertibility result \cite{SP09}.

\section{Minimal Sequences}
\setcounter{equation}{0}
The {\bf symbolic dynamical system} $(\Omega,\sigma),$ where $\Omega
:= \{0,1\}^{\ZZ}$ has the product topology and $\sigma:\Omega
\rightarrow \Omega$ is the shift homeomorphism $\sigma(b)(j) =
b(j-1), \, b \in \Omega,$ belongs to the class of dynamical systems
introduced by Bebutov in \cite{BE1940}. Its subsystems $(X,\sigma)$
correspond to nonempty closed invariant $X \subseteq \Omega.$
Elements in $\Omega$ are binary sequences and the sets $U_m(b) := \{
\, a \in \Omega \, : \, a(k) = b(k), -m < k < m \, \},$ $b \in
\Omega, \, m \in \NN$ are a basis for the product topology. Orbits
$O(b) := \{ \, \sigma^k(b) : k \in \ZZ \, \}$ are (shift) invariant
and orbit closures $\overline O(b)$ are closed and invariant.
\begin{lem}
\label{lem2.1} If $B(S,\Lambda)$ is a Riesz set and if $b$ is a
nonzero sequence in $\overline O(\chi_{_{\Lambda}})$ then
$B(S,supp(b))$ is a Riesz set.
\end{lem}
\begin{pf}
Fix $\epsilon_1 > 0.$ Then $B(S,\Lambda)$ satisfies the inequality
in condition (1) of Lemma (\ref{lem1.1}) iff $B(S,\Lambda_f)$
satisfies this inequality for every finite $\Lambda_f \subseteq
\Lambda.$ The result then follows from the definition of orbit
closure and product topology.
\end{pf}

\medskip

\noindent A nonempty closed invariant $X \subset \Omega$ is called a
minimal set if it is minimal with respect to these properties.
Zorn's lemma ensures that every nonempty closed invariant set
contains a minimal set. If $X$ is a minimal set and $b \in X$ then
$\overline O(b) = X.$ A {\bf minimal sequence} is a binary sequence
$b$ such that $\overline O(b)$ is a minimal set.
\begin{defin}
\label{defin2.1} $\Lambda \subset \ZZ$ is {\bf syndetic} if there
exists $n \in \NN$ such that $\Lambda +  F_n = \ZZ,$ thick if for
every $n \in \NN$ there exists $k \in \ZZ$ such that $k + F_n
\subset \Lambda,$ and piecewise syndetic if $\Lambda = \Lambda_s
\cap \Lambda_t$ where $\Lambda_s$ is syndetic and $\Lambda_t$ is thick.
\end{defin}
\begin{lem}
\label{lem2.2}
If $\ZZ = \bigcup_{\, i=1}^{\, n} \Lambda_i$ then one of
the $\Lambda_i$ is piecewise syndetic.
\end{lem}
\begin{pf}
Theorem 1.23 in \cite{FU81}.
\end{pf}
\begin{lem}
\label{lem2.3} If $\Lambda_p$ is piecewise syndetic then there
exists a syndetic set $\Lambda$ such that
$\chi_{_{\Lambda}} \in \overline O(\chi_{_{\Lambda_p}}).$
\end{lem}
\begin{pf}
$\Lambda_p = \Lambda \cap \Lambda_t$ where $\Lambda$ is syndetic and
$\Lambda_t$ is thick. Then
$\chi_{_\Lambda} \in \overline O(\chi_{_{\Lambda_p}})$
follows from the definitions of thick sets, orbit closures, and product topology.
\end{pf}
\begin{cor}
\label{cor2.1}
For every $S \in \cL^{+}(\TT)$ the following conditions are equivalent:
\begin{enumerate}
\item $B(S,\ZZ)$ is a finite union of Riesz sequences,
\item there exists a {\bf syndetic} subset $\Lambda \subset \ZZ$
such that $B(S,\Lambda)$ is a Riesz sequence.
\end{enumerate}
\end{cor}
\begin{pf}
(2) implies (1): If $\Lambda$ is syndetic there exists $n \in \NN$
with $\Lambda + F_n = \ZZ.$ Then $B(S,\ZZ)$ is the union of the
Riesz sequences $B(S,\Lambda + k), \, k \in F_n.$
\\
\noindent (1) implies (2): If $B(S,\ZZ)$ is a finite union of Riesz
sets then Lemma (\ref{lem2.2}) implies that there exists a piecewise
syndetic $\Lambda_p$ such that $B(S,\Lambda_p)$ is a Riesz set. Then
Lemma (\ref{lem2.3}) implies there exists a syndetic $\Lambda$ such
that $\chi_{_{\Lambda}} \in \overline O(\Lambda_p).$ Since $\Lambda
= supp(\chi_{_\Lambda}),$ Lemma (\ref{lem2.1}) implies that
$B(S,\Lambda)$ is a Riesz sequence.
\end{pf}
\begin{lem}
\label{lem2.4} If $\Lambda$ is syndetic and $b \in \overline
O(\chi_{_\Lambda})$ then $supp(b)$ is syndetic.
\end{lem}
\begin{pf}
Since $\Lambda$ is syndetic there exists $n \in \NN$ with $\Lambda +
F_n = \ZZ.$ Therefore
$$supp \left( \sigma^{k}(\chi_{_\Lambda}) \right) + F_n = supp \left(
\chi_{_\Lambda} \right) + k + F_n = \ZZ, \ \ k \in \ZZ,$$
so the definition of orbit closure implies $supp(b) + F_n = \ZZ$
whenever $b \in \overline O(\chi_{_\Lambda}).$
\end{pf}
For $b \in \Omega$ define the function $\theta_b : \ZZ \rightarrow
\Omega$ by
$\theta_b(k) = \sigma^k(b), \ \ k \in \ZZ.$
\begin{defin}
\label{defin2.2} $b \in \Omega$ is almost periodic if
$\theta_{b}^{-1}(U_m(b))$ is syndetic for every $m \in \NN.$
\end{defin}
\begin{lem}
\label{lem2.5} A sequence is minimal iff it is almost periodic.
\end{lem}
\begin{pf}
Gottschalk and Hedlund proved this in \cite{GH55}, Theorems (4.05)
and (4.07).
\end{pf}
\begin{cor}
\label{cor2.2}
If $b$ is a nonzero minimal sequence then $supp(b)$ is syndetic.
\end{cor}
\begin{pf}
Choose $k \in supp(b)$ and set $m = |k| + 1.$ Lemma (\ref{lem2.5})
implies that $b$ is almost periodic therefore there exists $n \in
\NN$ such that $\theta_{b}^{-1}(U_m(b)) + F_n = \ZZ.$ Therefore
$supp(b) + F_n = \ZZ$ so $supp(b)$ is syndetic.
\end{pf}
{\bf Proof of Theorem \ref{thm_main}}
(1) equivalent to (2): This follows from Corollary (\ref{cor2.1}).
\\
(3) implies (2): Since $\Lambda$ is nonempty $\chi_{_\Lambda}$ is a nonzero minimal
sequence and hence Corollary (\ref{cor2.2}) implies that $\Lambda = supp(\chi_{_\Lambda})$ is syndetic.
\\
(2) implies (3): Zorn's lemma implies that there exists a
minimal set $X \subseteq \overline O(\chi_{_\Lambda}).$ Then choose
$b \in X.$ Then $b$ is a minimal sequence. Lemma (\ref{lem2.1}) implies that
$B(S,supp(b))$ is a Riesz sequence. Lemma (\ref{lem2.4}) implies that $supp(b)$
is syndetic and hence $supp(b)$ is nonemepty. Then (3) follows from the fact
that $\chi_{_{supp(b)}} = b.$
\begin{defin}
\label{defin2.3} A subset $\Lambda \subset \ZZ$ is a Bohr set if
there exists a compact abelian group $G,$ a homomorphism $\psi : \ZZ
\rightarrow G$ with $\overline {\psi(\ZZ)} = G,$ and a nonempty open
subset $U \subset G$ such that $\Lambda = \psi^{-1}(U).$
\end{defin}
If $\Lambda$ is a Bohr set then $\chi_{_\Lambda}$ is a nonzero
minimal sequence. These sets generalize sets having the form $n\ZZ +
m$ where $n \in \NN$ and $m \in \ZZ$ and are unions of the Bohr sets
defined by Ruzsa (\cite{GR09}, Definition 2.5.1) who studied their
number theoretic properties. They are named after Harald Bohr, who
pioneered the theory of (uniformly) almost periodic functions
\cite{BO52}, and are related to the Bohr compactification used
by Dutkay, Han, and Jorgensen in their study of spectral pairs
\cite{DU09}. The following extension of Result 2 utilizes
spectral properties of Bohr sets.
\begin{thm}
\label{thm2} If $S$ is a Cantor set with $\mu(S) > 0$ and $\Lambda$
is a Bohr set then $B(S,\Lambda)$ is not a Riesz set.
\end{thm}
\begin{pf}
Without loss of generality we can assume that $\Lambda = \psi^{-1}(U)$ where
$U \subseteq G$ is an open set that contains $0 \in G$ and choose an open
subset $V \subseteq G$ that contains $0$ and satisfies $V - V \subseteq U.$
Set
$f := \chi_{_V} * \chi_{_{-V}}$
and
$g := f \circ \psi \in \ell^{\infty}(\ZZ).$
Then $supp(g) \subset \Lambda$ and
$g$ equals the Fourier transform $\widehat \nu$ of the
positive measure $\nu$ on $\TT$  given by
\begin{equation}
\label{eqnnu}
\nu = \sum_{\gamma \in \widehat G} \widehat f(\gamma) \,
\delta_{_{\gamma(\psi(1))}}, \ \ \ \ \widehat f(\gamma) = | \widehat
\chi_{_V}(\gamma)|^2
\end{equation}
where $\widehat G$ is the Pontryagin dual of $G$ and $\widehat f \in
\ell^{2}(\widehat G)$ is the Fourier transform of $f.$ Let $\epsilon
> 0.$ It suffices to construct $h \in L^2(\TT)$ such that
$||P_S \, (\nu*h) || < \epsilon \, || \nu*h ||$
since $P_{_\Lambda} (\nu*h) = \nu*h.$ Partition $\widehat G =
\Gamma_1 \cup \Gamma_2$ where $\Gamma_1$ is finite,
let $\nu_i$ be the component of $\nu$ supported on $\Gamma_i, i = 1,2,$ and
set $\alpha := \sum_{\gamma \in \Gamma_2} \widehat f(\gamma)$ and
$\beta := \sum_{\gamma \in \Gamma_1} \widehat f(\gamma)^2.$ Since
$S$ is nowhere dense $supp(\nu_1) + S \neq \TT$ so there exists
$h \in L^2(\TT)$ such that $||h|| = 1$ and
$supp(h)$ is contained in an arc $I \subset \TT$ that is disjoint
from $supp(\nu_1) + S$ and such that the intervals $I + \gamma,
\gamma \in \Gamma_1$ are mutually disjoint.
Then
$||P_S (\nu * h)|| = ||P_S (\nu_2 * h)|| \leq || \nu_2 * h || \leq
\alpha,$
and
$|| \nu*h || \geq || \nu_1*h || = \beta.$
The result follows by choosing $\Gamma_1$ so $\alpha < \epsilon
\beta$ which is possible since as $\Gamma_1$ increases $\alpha
\rightarrow 0$ and $\beta \rightarrow f(1) > 0.$
\end{pf}
Bohr minimal sequences are simple. We discuss methods to construct
more sophisticated minimal sequences. For nonempty invariant $X, Y
\subseteq \Omega,$ a function $\zeta : X \rightarrow Y$ is
equivariant if $\zeta \circ \sigma = \sigma \circ \zeta.$ For $m \in
\NN$ every function $c : \{0,1\}^{\{-m+1,...,m-1\}} \rightarrow
\{0,1\}$ defines the function $\zeta_c : \Omega \rightarrow \Omega$
by
\begin{equation}
\label{zeta}
    \zeta_c(b)(k) = c \left( R_{\{-m+1,...,m-1\}} (\sigma^k(b)) \right),
    \ \ b \in \Omega, \, k \in \ZZ.
\end{equation}
Furthermore, for every nonempty closed invariant $X \subseteq
\Omega$ the restriction $\zeta_c : X \rightarrow \Omega$ is
continuous and equivariant and every continuous equivariant $\zeta :
X \rightarrow \Omega$ equals $\zeta_c$ for some $c.$ Equivariant
images of minimal sets and sequences are minimal.

\medskip

\noindent The Thue-Morse minimal sequence
$b = \cdots 10010110.0110100110010110 \cdots,$
introduced in \cite{TH06}, \cite{MO21}, can be constructed
using substitutions $0 \rightarrow 01$ and $1 \rightarrow 10.$ Its
orbit closure $X = \overline O(b)$ admits a unique invariant
ergodic probability measure $\lambda$ \cite{KE68}.
The spectrum of the unitary operator
$(U_{\sigma}f)(x) = f(\sigma(x)), \, f \in L^2(X,\lambda)$ admits
a Reisz product representation, has no point components, and is
supported on a dense set of measure zero \cite{MA26}, \cite{Q87}.

\section{Research Directions}
\setcounter{equation}{0}

We suggest three questions, related to the material in
this paper, as directions towards a solution of the FCE.
In this section we assume that $S \in \cL^{+}(\TT)$ is a Cantor
set such that $\chi_{_S} \not \in W_{2,2}^{\tau}$ for all $\tau > 0$
and that $\chi_{_\Lambda}$ is a nonzero minimal sequence.
We let $\cM(\Lambda),$ $\cP(\Lambda)$ denote the set of measures, pseudomeasures, respectively,
on $\TT$ whose Fourier transforms are supported on $\Lambda,$ see (\cite{LA71},4.2).

\medskip

\noindent {\bf Question 1} What properties of a pair $(S,\Lambda),$
determine whether or not $B(S,\Lambda)$ is a Riesz sequence? Such
properties include the rate of decay of the restriction of
$\widehat \chi_{_S}$ to $\Lambda,$ and the sumsets $S + supp(\nu)$
where $\nu \in \cM(\Lambda)$ or $\nu \in \cP(\Lambda).$
Of particular interest are pairs where $\chi_{_\Lambda}$ is a substitution minimal
sequence because their spectral properties have been intensively studied \cite{Q87}.

\medskip

\noindent {\bf Question 2} What is the spectrum of $P_{_S} + P_{_\Lambda}?$
Condition (2) in Lemma (\ref{lem1.1}) implies that $B(S,\Lambda)$ is
a Riesz sequence iff this spectrum is bounded below by a
positive number. Let $\cA(P_{_S},P_{\Lambda})$ denote the $C^*$-subalgebra
of $\cB(L^2(\TT))$ generated by $P_{_S}$ and $P_{_\Lambda}$. A
standard result (\cite{FI96}, 8.5.5) shows that $\cA(P_{_S},P_{\Lambda})$
equals a homomorphic image of a specific crossed-product $C^*$-algebra and implies
that $\cA(P_{_S},P_{\Lambda})$ is determined by the spectrum of $P_{_S} + P_{_\Lambda}.$

\medskip

\noindent {\bf Question 3} What are the spectrums of submatrices of
the Laurent operators
\\
$L_{_S} : \ell^2(\ZZ) \rightarrow \ell^2(\ZZ)$ defined
by $L_{_S} f = \widehat \chi_{_S} * f?$ Spielman and Srivastava's
algorithm \cite{SP09} may provide an efficient method to compute these spectrums.
Of particular interest are Cantor sets having the form $S = \bigcap_{n \in \NN} S_n$
where each $S_n$ is obtained by deleting a large number of equally spaced, equal length
open arcs from $\TT.$ This construction was suggested to the author by Alexander Olevskii
as a method of constructing Cantor sets $S$ such that $\widehat \chi_{_S}$ decays slowly and
is easily computable.

\medskip

\noindent {\it Acknowledgments} We thank Eric Weber for introducing
the beguiling Kadison-Singer problem and Peter Casazza, Ilya Chrishtal,
Ole Christiansen, Palle Jorgensen, Anders Mauritzen, and
Alexander Olevskii for their advice and encouragement.

\end{document}